\documentclass[12pt]{article}
\usepackage{amssymb,amsmath,amsfonts,amsthm}
\usepackage{bm}
\usepackage{cite}
\usepackage{url}

\def\e{\mathrm e}                                               
\def\om{{\mathrm r}}                                            
\def\RR{\mathop{\cal R}\nolimits}                               
\def\cdc{,\ldots,}                                              
\def\1n{1\cdc n}                                                
\def\eq#1{\begin{equation}#1\end{equation}}                     
\def\eqs*#1{\begin{eqnarray*}#1\end{eqnarray*}}                 
\def\eqss#1{\begin{eqnarray}#1\end{eqnarray}}                   
\newtheorem{thm}{Theorem}{\bfseries}{\itshape}                  
\newtheorem{corol}{Corollary}{\bfseries}{\itshape}              
\newtheorem{defin}{Definition}{\bfseries}{\upshape}             
\newtheorem{lemma}{Lemma}{\bfseries}{\itshape}                  
\def\proof{{\noindent\bf Proof. }}                              
\def\R{{\mathbb R}}                                             
\def\T{{\xz\rm\scriptscriptstyle T}\xz}                         
\def\LW{{\rm\scriptscriptstyle LW}}                             
\def\xy{\hspace{.07em}}                                         
\def\xz{\hspace{-.07em}}                                        
\def\hy{\hspace{.04em}}                                         
\def\hz{\hspace{-.04em}}                                        
\def\ms{\mathstrut}                                             
\def\diag{\operatorname{diag}}                                  
\def\_#1{{^{}_{#1}}}                                            
\def\jj{{\bar\jmath\bar\jmath}}                                 
\def\ii{{\bar\imath\hspace{.075ex}\bar\imath}}                  
\def\uv{{\hy\overline{\hz uv\hz}\hy}}                           
\def\vu{{\hy\overline{\hz vu\hz}\hy}}                           
\def\vv{{\hy\overline{\hz vv\hz}\hy}}                           
\def\u{{\hy\overline{\hz u\hz}\hy}}                             
\def\v{{\hy\overline{\hz v\hz}\hy}}                             
\def\l{\ell}                                                    
\DeclareMathSymbol{\ell}{\mathord}{letters}{96}                 
\def\dr{d^{\hspace*{.08em}\om}}                                 
\def\dLW{d^{\hspace{.05em}\LW}\xz}                              
\def\tdLW{\tilde d^{\hspace{.05em}\LW}\xz}                      
\def\Up#1{\vspace{-#1em}}                                       
\def\L{{\cal L}}                                                
\def\Ker{\operatorname{Ker}}                                    
\def\rank{\operatorname{rank}}                                  
\def\too{\!\to\!}                                               

\author{Pavel Chebotarev\footnotemark[1]
        \and  R.B. Bapat\footnotemark[2]
        \and   R. Balaji\footnotemark[3]
}

\title{\vspace{-.0em}Simple expressions for the long walk distance}
\date{}

\def\thefootnote{\fnsymbol{footnote}}
\begin{document}
\def\baselinestretch{0.94}

\footnotetext[1]{{Institute of Control Sciences of the Russian Academy of Sciences,}
                               {65 Profsoyuznaya Street, Moscow 117997, Russia;}
                               e-mail: {\footnotesize\tt chv@member.ams.org.}\\
                               The work of the author has been supported by the RAS Program ``Network Control under Uncertainty''.}
\footnotetext[2]{{Indian Statistical Institute, New Delhi 110016, India;}
                               e-mail: {\footnotesize\tt rbb@isid.ac.in.}\\
                               The work of the author has been supported by the JC Bose Fellowship, Department of Science and Technology, Government of India.}
\footnotetext[3]{{Department of Mathematics, Indian Institute of Technology-Madras, Chennai-36, India;}
                               e-mail: {\footnotesize\tt balaji5@iitm.ac.in.}}
\maketitle
\def\thefootnote{\arabic{footnote}}
\vspace{-2em}
\begin{abstract}
The walk distances in graphs are defined as the result of appropriate transformations of the $\sum_{k=0}^\infty(tA)^k$ proximity measures, where $A$ is the weighted adjacency matrix of a connected weighted graph and $t$ is a sufficiently small positive parameter. The walk distances are graph-geodetic, moreover, they converge to the shortest path distance and to the so-called long walk distance as the parameter $t$ approaches its limiting values. In this paper, simple expressions for the long walk distance are obtained. They involve the generalized inverse, minors, and inverses of submatrices of the symmetric irreducible singular M-matrix $\L=\rho I-A,$ where $\rho$ is the Perron root of~$A.$


\end{abstract}

\section{Introduction}

The walk distances for graph vertices are a parametric family of graph distances proposed in~\cite{Che11AAM}.
Along with their modifications they generalize~\cite{Che12DAM} the logarithmic forest distances \cite{Che11DAM}, resistance distance, shortest path distance, and the weighted shortest path distance.
The walk distances are graph-geodetic: for a distance\footnote{In this paper, a \emph{distance\/} is assumed to satisfy the axioms of metric.} $d(i,j)$ on the vertices of a graph $G$\/ this means that $d(i,j)+d(j,k)=d(i,k)$ if and only if every path in $G$ connecting $i$ and $k$ visits~$j.$

The long walk distance, $\dLW(i,j)$, is obtained by letting the parameter of the walk distances go to one of its limiting values (approaching the other limiting value leads to the shortest path distance). A number of expressions for $\dLW(i,j)$ are given in~\cite{Che12DAM}. In this paper, we obtain simple expressions in terms of the matrix $\L=\rho I-A,$ where $A$ is the weighted adjacency matrix of a connected weighted graph $G$ on $n$ vertices and $\rho$ is the Perron root of~$A.$ $\L$~is a symmetric irreducible singular M-matrix, so $\rank\L=n-1$ and $\L$ is positive semidefinite. In \cite{Che12DAM}, $\L$ was called the \emph{para-Laplacian matrix\/} of~$G.$ The expressions presented in this paper generalize some classical expressions for the resistance distance (cf.~\cite{BapatSivasubramanian11}).
They enable one to conclude that the long walk distance can be considered as the counterpart of the resistance distance obtained by replacing the Laplacian matrix $L\!=\!\diag(A\bm1)\!-\!A$ and the vector~$\bm1$ (of $n$ ones) which spans $\Ker L$ with the para-Laplacian matrix $\L$ and the eigenvector of $A$ spanning $\Ker\L.$ If $A$ has constant row sums, then $\L=L$ and these distances coincide.

\section{Notation}
\label{s_notat}

In the graph definitions we mainly follow~\cite{Harary69}.
Let $G$ be a weighted multigraph (a weighted graph where multiple edges are allowed) with vertex set $V(G)=V,$ $|V|=n>1,$
and edge set~$E(G)$. Loops are allowed; throughout the paper we assume that $G$ is connected. For brevity, we will call $G$ a \emph{graph}.
For ${i,j\in V,}$ let $n_{ij}\in\{0,1,\ldots\}$ be the number of edges incident to both $i$ and $j$ in~$G$; for every ${q\in\{\1n_{ij}\}}$, $w_{ij}^q>0$ is the weight of the $q\/$th edge of this type. Let
\eqs*{
a_{ij}=\sum_{q=1}^{n_{ij}}w_{ij}^q
}
(if $n_{ij}=0,$ we set $a_{ij}=0$) and $A=(a_{ij})_{n\times n};$ $A$~is the symmetric \emph{weighted adjacency matrix\/} of~$G$.
In this paper, all matrix entries are indexed by the vertices of~$G.$

The \emph{weighted Laplacian matrix\/} of $G$ is $L\!=\!\diag(A\bm1)\!-\!A,$ where $\bm1$ is the vector of $n$ ones.

\smallskip
For $v_0,v_m\in V,$ a $v_0\to v_m$ \emph{path\/} (\emph{simple path}) in $G$ is an alternating sequence of vertices and edges $v_0,\e_1,v_1\cdc\e_m,v_m$ where all vertices are distinct and each $\e_i$ is an edge incident to $v_{i-1}$ and~$v_i.$ The unique $v_0\to v_0$ path is the ``sequence''\,$v_0$ having no edges.

Similarly, a $v_0\to v_m$ \emph{walk\/} in $G$ is an \emph{arbitrary\/} alternating sequence of vertices and edges $v_0,\e_1,v_1\cdc\e_m,v_m$ where each $\e_i$ is incident to $v_{i-1}$ and~$v_i.$
The \emph{length\/} of a walk is the number $m$ of its edges (including loops and repeated edges). The \emph{weight\/} of a walk is the product of the $m$ weights of its edges. The weight of a set of walks is the sum of the weights of its elements. By definition, for any vertex $v_0$, there is one $v_0\to v_0$ walk $v_0$ with length $0$ and weight~1.

Let $r_{ij}$ be the weight of the set $\RR^{ij}$ of all $i\to j$ walks in $G$, provided that this weight is finite. $R=R(G)=(r_{ij})_{n\times n}$ will be referred to as the \emph{matrix of the walk weights}.

By $\dr(i,j)$ we denote the \emph{resistance distance\/} between $i$ and $j$, i.e., the effective resistance between $i$ and $j$ in the resistive network whose line conductances equal the edge weights $w_{ij}^q$ in~$G$. The resistance distance has several expressions via the generalized inverse, minors, and inverses of the submatrices of the weighted Laplacian matrix of $G$ (see, e.g.,~\cite{KleinRandic93,Bapat99RD} or the papers by Subak-Sharpe and Styan published in the 60s and cited in~\cite{Che12DAM}).

\begin{defin}
\label{d_g-d}
{\rm
For a multigraph $G$ with vertex set $V,$ a function $d\!:V\!\times\!V\to\R$ is called \emph{graph-geodetic\/} provided that
$d(i,j)+d(j,k)=d(i,k)$ holds if and only if every path in $G$ connecting $i$ and~$k$ contains~$j$.
}
\end{defin}

Graph-geodetic functions can also be called \emph{bottleneck additive\/} or \emph{cutpoint additive}. 

\def\baselinestretch{1.0}
\Up{.6}\section{The walk distances}
\label{s_prel}

For any $t>0,$ consider the graph $tG$ obtained from $G$ by multiplying all edge weights by~$t.$
If the matrix of the walk weights of $tG,\,$ $R_t=R(tG)=(r_{ij}(t))_{n\times n},$ exists, then\footnote{In the more general case of weighted \emph{digraphs}, the $ij$-entry of the matrix $R_t-I$ is called the \emph{Katz similarity\/} between vertices $i$ and~$j$. Katz~\cite{Katz53} proposed it to evaluate the social status taking into account all $i\too j$ paths.
}
\eqs*{
R_t=\sum_{k=0}^\infty(tA)^k=(I-tA)^{-1},
}
where $I$ denotes the identity matrix of appropriate dimension.

By assumption, $G$ is connected, while its edge weights are positive, so $R_t$ is positive whenever it exists.
The walk distances can be introduced as follows.

\begin{defin}
\label{d_walkD}
{\rm
For a connected graph $G,$ the \emph{walk distances\/} on~$V(G)$ are the functions $d_t(i,j)\!:V(G)\!\times\! V(G)\to\R$ and the functions
positively proportional to them$,$ where
\eqss{
\label{e_d1}
d_t(i,j)
=-\ln\biggl(\xz\frac{r_{ij}(t)}{\sqrt{r_{ii}(t)\,r_{\!jj}(t)}}\xz\biggr).
}
}
\end{defin}

\begin{lemma}[\!\!\cite{Che11AAM}]
\label{l_di}
For any connected $G,$ if ${R_t=(r_{ij}(t))}$ exists$,$ then \eqref{e_d1} determines a graph-geodetic distance on~$V(G).$
\end{lemma}

Regarding the existence (finiteness) of $R_t,$ since $G$ is connected$,$ $A$ is irreducible, so the Perron-Frobenius theory of nonnegative matrices provides the following result.
\begin{lemma}
\label{l_finite}
For any weighted adjacency matrix $A$ of a connected graph $G,$ the series
$R_t=\sum_{k=0}^\infty(tA)^k$ with $t>0$ converges to $(I-tA)^{-1}$ if and only if\/ $t<\rho^{-1},$ where $\rho=\rho(A)$ is the spectral radius of\/~$A.$ Moreover$,$ $\rho$ is an eigenvalue of~$A;$ as such $\rho$ has multiplicity~$1$ and a positive eigenvector.
\end{lemma}

Thus, the walk distance $d_t(i,j)$ with parameter $t$ exists if and only if $0<t<\rho^{-1}.$

A topological interpretation of the walk distances in terms of the weights of routes and circuits in $G$ is given in~\cite{CheDeza11}.

\section{The long walk distance}

The \emph{long walk distance\/} is defined \cite{Che12DAM} as the limit of the walk distance \eqref{e_d1} as $t\to(\rho^{-1})^-$ with a scaling factor depending on $t,$ $\rho,$ and~$n.$ Namely,
\eqs*{
\dLW(i,j)&\stackrel{\rm def}
 =&\lim_{t\to(\rho^{-1})^-}\frac{2d_t(i,j)}                                 {n\rho^2(\rho^{-1}-t)}\\
&=&\lim_{t\to(\rho^{-1})^-}\frac{\ln r_{ii}(t)+\ln r_{jj}(t)-2\ln r_{ij}(t)}{n\rho^2(\rho^{-1}-t)},\quad i,j\in V.
}

It is shown \cite{Che12DAM} that this limit always exists and defines a graph-geodetic metric.

There is a two-way connection between the long walk distance and the resistance distance. In particular, the following theorem holds.

\begin{thm}[\!\!\cite{Che12DAM}]
\label{t_rel}
Let $G$ be a connected graph on $n$ vertices with weighted adjacency matrix $A;$ let $p=(p_1\cdc p_n)^\T$ be the Perron vector of~$A.$ Suppose that $p'=(p'_1\cdc p'_n)^\T=\sqrt{n}p/\|p\|_2,$ where $\|p\|_2=(\sum_{i=1}^np_i^2)^{1/2},$ and $P'=\diag p'.$ Then the long walk distance in $G$ coincides with the resistance distance in the graph $G'$ whose weighted adjacency matrix is~$P'AP'.$
\end{thm}

Using Theorem\:\ref{t_rel} and two classical expressions for the resistance distance one has:

\begin{corol}[\!\!\cite{Che12DAM}]
\label{c_frResi}
In the notation of Theorem\/\:{\rm\ref{t_rel},}
\eqs*{
\dLW(i,j)
&=&\frac{\det\xy\xy(L'_\ii)_\jj}{\det L'_\vv},\quad j\ne i,\quad\forall\, v\in V,\\
\dLW(i,j)
&=&\l'^{\bm-}_{ii}+\l'^{\bm-}_{jj}-2\l'^{\bm-}_{ij},
}
where $L'=P'\L\xy P',$ $\L=\rho I-A,\,$ $X_\uv$ is $X$ with row $u$ and column $v$ removed\/$,$
and\/ $L'^{\bm-}=(\l'^{\bm-}_{ij})$ is any g-inverse\/\footnote{$Z$ is a \emph{g-inverse\/} of $X$ whenever $X=XZX.$ A simple choice of $L'^{\bm-}$ is $(L'+\frac1n\bm{11}^\T)^{-1}$ \cite[Theorem\:10.1.4]{RaoMitra71}.} of\:$L'.$
\end{corol}

Note that by Lemma\:13 in \cite{Che12DAM}, $L'$ is the weighted Laplacian matrix of the connected graph $G'$ introduced in Theorem\:\ref{t_rel}. Therefore, the Matrix tree theorem implies that $\det L'_\vv,$ the denominator of the above determinant formula for $\dLW(i,j),$ is strictly positive.

Now we obtain simple expressions for $\dLW(i,j)$ in terms of the matrix~$\L=\rho I-A.$

\section{Simple expressions for the long walk distance}

\begin{thm}
\label{t_LW_I-A}
In the notation of Theorem\/\;$\ref{t_rel}$ and
Corollary\/\;{\rm\ref{c_frResi},}
for all\/ $i,j\in V$ such that $j\ne i,$
\eqss{
\label{e_dLW_L1}
\dLW(i,j)
&=&\frac{\det\xy\xy(\L_{\xy\ii})_\jj}{p'^2_j\det \L_{\xy\ii_{\ms}}},\\
\label{e_dLW_L2}
\dLW(i,j)
&=&z^\T(i,j)\xy\L^{\bm-}z(i,j)^{\ms},\\
\label{e_dLW_L3}
\dLW(i,j)
&=&z_{\xy\u}^\T\xy(i,j)\hy(\L_\vu)^{-1}z\_\v(i,j)^{\ms},\quad \forall\, u,v\in V,
}
where
$\L^{\bm-}$ is any g-inverse of~$\L$ and\/ $z(i,j)$ is the $n$-vector whose $i$th element is $1/p'_i,$ $j$th element is\, $-1/p'_j,$ and the other elements are\/~$0$.
\end{thm}

\proof
By virtue of Corollary\:\ref{c_frResi}, for $j\ne i$ we have
\eqs*{
\dLW(i,j)
=\frac{\det\xy((P'\L P')_\ii)_\jj}{\det\xy(P'\L P')_\ii}
=
\frac{(\prod_{k\not\in\{i,\xy j\}}p'_k)^2\det\xy(\L_\ii)_\jj}
     {(\prod_{k\ne i}             p'_k)^2\det\xy \L_\ii}
=\frac{\det\xy\xy(\L_{\xy\ii})_\jj}{p'^2_j\det \L_{\xy\ii_{\ms}}},
}
that is, \eqref{e_dLW_L1} holds. To prove~\eqref{e_dLW_L2}, consider the matrix\footnote{$H$ generalizes the \emph{zero-axial matrix\/} $Z$ studied in~\cite{StyanSubak-Sharpe97}.} $H=(h_{ij})\in\R^{n\times n}$ such that
\eq{
\label{e_uij}
h_{ij}=
\begin{cases}
-\dfrac{ p'_i\det\xy\xy(\L_{\xy\ii})_\jj}
       {2p'_j\det\L_{\xy\ii_{\ms}}\phantom{(_\jj)}},& j\ne i,\\
0,                                                  & j=i.
\end{cases}
}

\begin{lemma}
\label{l_BaBa}
$H$ is a g-inverse of $\L$\,$:\,$ $\L\xy H\L=\L$.
\end{lemma}

\noindent\textbf{Proof.}
Let $D_{G'}^{\xy\om}$ be the matrix of resistance distances in the graph $G'$ introduced in Theorem\:\ref{t_rel}. Then, by \cite[Theorem\:10.1.4]{RaoMitra71},
\eq{
\label{e_DrG'}
D_{G'}^{\xy\om}=\l^{\bm-}_{G'}\bm1^\T+\bm1(\l^{\bm-}_{G'})^\T-2L^{\bm-}_{G'},
}
where $L^{\bm-}_{G'}$ is any g-inverse of $L_{G'},$ $L_{G'}$ is the weighted Laplacian matrix of $G',$ and $\l^{\bm-}_{G'}$ is the column vector containing the diagonal entries of $L^{\bm-}_{G'}.$ As $L_{G'}\bm1=\bm0$ and $\bm1^\T L_{G'}=\bm0^{\xy\T},$ \eqref{e_DrG'} implies that
\eqs*{
L_{G'}(-\tfrac12D_{G'}^{\xy\om})L_{G'}=L_{G'}.
}
Since by Lemma\:13 in \cite{Che12DAM}, $L_{G'}=P'\L P'$ and by Theorem\:\ref{t_rel}, $D_{G'}^{\xy\om}=D_G^{\xy\LW}$ (where $D_G^{\xy\LW}$ is the matrix of long walk distances in~$G$), we have
\eq{
\label{e_mnog1}
P'\L P'(-\tfrac12D_G^{\xy\LW})P'\L P'=P'\L P'.
}

It follows from~\eqref{e_dLW_L1} and~\eqref{e_uij} that
\eqs*{
H
=P'(-\tfrac12D_G^{\xy\LW})P',
}
moreover, $P'$ is invertible. Thus, \eqref{e_mnog1} yields
$\L\xy H\L=\L,$ as desired.
\qed

\smallskip
By~\eqref{e_dLW_L1} and~\eqref{e_uij},
\eq{
\label{e_du}
d^{\xy\LW}(i,j)
=-2\frac{h_{ij}}{p'_i\xy p'_j}
=z^\T(i,j)\, H\xy z(i,j),\quad j\ne i,
}
where $z(i,j)$ is the $n$-vector whose $i$th element is $1/p'_i,$ $j$th element is\, $-1/p'_j,$ and the other elements are\/~$0$.

Let $i$ and $j\ne i$ be fixed and let $z=z(i,j).$
Lemma\:\ref{l_finite} implies that $0$ is a simple eigenvalue of $\L,$ while $p$ spans $\Ker\L.$ Hence, $\rank\L=n-1$ and the range of $\L$ is~$p^\bot,$ the orthogonal compliment of~$p.$ Since $z\in p^\bot,$ it follows that $z=\L\xy x$ for some $n$-vector~$x$. Consequently, if $\L^{\bm-}$ is a g-inverse of $\L,$ then
\eqs*{
 z^\T(i,j)\xy\L^{\bm-}z(i,j)
=x^\T\L      \L^{\bm-}\L\xy x
=x^\T                 \L\xy x
=z^\T(i,j)\xy H       z(i,j)
=d^{\xy\LW}(i,j),
}
where the last two transitions follow from Lemma\:\ref{l_BaBa} and Eq.\:\eqref{e_du}. This proves~\eqref{e_dLW_L2}.

Since $\L$ is a singular M-matrix, the principal minors of $\L$ are strictly positive. By \cite[Lemma\:7]{Che12DAM} (which follows from Theorem\:3.1 in~\cite{Meyer74}) $l^{vu}/l^{vv}=p_u/p_v,$ where $l^{vu}$ is the cofactor of the entry $l_{vu}$ of~$\L.$ Hence, $\L_\vu$ is non-singular.
In view of \eqref{e_dLW_L2}, proving \eqref{e_dLW_L3} reduces to showing that
\eq{
\label{e_LLLL}
\L\L_{(uv)}^{\bm\sim}\L=\L,
}
where $\L_{(uv)}^{\bm\sim}\!\in\!\R^{n\times n}$ is the matrix that has all zeros in row $u$ and column~$v$ and satisfies $\bigl(\L_{(uv)}^{\bm\sim}\bigr)\_\uv=(\L_\vu)^{-1}.$ Eq.\:\eqref{e_LLLL} can be derived from the general result \cite[Section\:11.2]{RaoMitra71}. Alternatively, the multiplication of block matrices verifies that for ${(l_{km})=\L}$ and $\L'=(l'_{km})\stackrel{\rm def}=\L\L_{(uv)}^{\bm\sim}\L$ it holds that $\L'_\uv=\L\_\uv,$ $l'_{u\v}=l_{u\v},$ $l'_{\u v}=l_{\u v},$ and $l'_{uv}=l_{u\v}(\L_\uv)^{-1} l_{\u v},$ where $l\_{u\v}$ and $l\_{\u v}$ are $\L$'s row $u$ and column $v$ with $l\_{uv}$ removed, respectively. Now to prove \eqref{e_LLLL} it remains to show that
\eq{
\label{e_idenL}
l_{u\v}(\L_\uv)^{-1} l_{\u v}=l\_{uv},
}
which is true due to the following lemma. 
\begin{lemma}
\label{l_iden1}
If $X=(x\_{km})_{n\times n}$ is a singular matrix$,$ while $\,X_\uv$ is non-singular$,$ then
\eqs*{
x_{uv}=x_{u\v}\,(X_{\uv})^{-1}\xy x_{\u v},
}
where $x_{ u\v}$ and\/ $x_{\u v}$ are $X\!$'s row $u$ and column $v$ with $x_{uv}$ removed\/$,$ respectively.
\end{lemma}

Lemma~\ref{l_iden1} is a special case of \cite[Eq.\,(6.2.1)]{Meyer00MAALA} and can alternatively be proved by first expanding $\det X$ along row $u$ and then expanding the cofactors of the entries in row $u$ (except for $x_{uv}$) along column~$v.$
Lemma~\ref{l_iden1} verifies \eqref{e_idenL} and thus, \eqref{e_LLLL} and~\eqref{e_dLW_L3}.
\qed

\bigskip
In conclusion, it is worth introducing a rescaled version of the long walk distance, 
$\tdLW(i,j)\!\stackrel{\rm def}=\!n\xy\|p\|_2^2\,\dLW(i,j),$ where, as before, $p$ is the Perron vector of~$A$ and $n$ is the number of vertices. It has the following desirable property: if for a sequence of connected graphs $(G^k)$ on $n$ vertices, we have\footnote{The superscript $k$ refers to the graph~$G^k.$} $a^k_{12}=a^k_{21}=1,$ $\,k=1,2\cdc$ whereas all other entries $a^k_{ij}$ tend to zero as $k$ goes to infinity, then $\tdLW(1,2)$ tends to~$1$ (while all other non-zero distances tend to infinity). If $G$ is balanced, i.e., all its vertices have the same weighted degree, then, obviously, both versions of the long walk distance coincide with the resistance distance: $\tdLW\equiv\dLW\equiv\dr.$


\end{document}